\documentclass[12pt]{amsart}
\usepackage{amssymb}

\newtheorem{thm}{Theorem}[section]

\newtheorem{lemma}[thm]{Lemma}
\newtheorem{cor}[thm]{Corollary}
\newtheorem{conj}[thm]{Conjecture}
\newtheorem{remark}[thm]{Remark}
\newenvironment{rem}{\begin{remark}\rm}{\end{remark}}
\newtheorem{example}[thm]{Example}
\newenvironment{ex}{\begin{example}\rm}{\end{example}}

\newcommand{\Edot}{{E_\bullet}}
\newcommand{\Fdot}{{F_\bullet}}

\newcommand{\Gdot}{{G_\bullet}}

\newcommand{\Fla}{{{\mathbb F}\ell_{\bf d}}}

\newcommand{\Span}[1]{{\langle #1 \rangle}}

\begin{document}

\title[Real enumerative geometry for flag manifolds]{Some real and unreal
enumerative geometry\\ for flag manifolds}   

\author{Frank Sottile}
\address{Department of Mathematics\\
        University of Wisconsin\\
        Van Vleck Hall\\
        480 Lincoln Drive\\
        Madison, Wisconsin 53706-1388\\
        USA}
\curraddr{Department of Mathematics\\
        University of Massachusetts\\
        Amherst, Massachusetts ??????\\
        USA}
\email{sottile@math.umass.edu}
\urladdr{http://www.math.umass.edu/\~{}sottile}
\date{17 July 2000}
\thanks{Research done in part while visiting IRMA in Strasbourg and
Universit\'e de Gen\`eve, and supported in part by Fonds National Suisse pour
la recherche} 
\thanks{Michigan Mathematics Journal, to appear.}
\dedicatory{To Bill Fulton on the occasion of his 60th birthday.}

\begin{abstract}
We present a general method for constructing real solutions to some problems
in enumerative geometry which gives lower bounds on the
maximum number of real solutions.
We apply this method to show that two new classes of enumerative
geometric problems on flag manifolds may have all their solutions be real
and modify this method to show that another class may have no real
solutions, which is a new phenomenon. 
This method originated in a numerical homotopy
continuation algorithm  adapted to the special Schubert calculus on
Grassmannians 
and in principle gives optimal numerical homotopy algorithms for finding
explicit solutions to these other enumerative problems.
\end{abstract}

\maketitle

\section*{Introduction}
For us, enumerative geometry is concerned with counting the geometric
figures of some kind that have specified position with respect to some
fixed, but general, figures.
For instance, how many lines in space are incident on four general (fixed)
lines?  (Answer: 2.) \ 
Of the figures having specified positions with respect to fixed
{\it real} figures, some will be real while the
rest occur in complex conjugate pairs, and the distribution between these
two types depends subtly upon the
configuration of the fixed figures.
Fulton~\cite{Fu96b} asked how many solutions to such a problem of enumerative 
geometry can be real and later with Pragacz~\cite{FuPr} reiterated this
question in the context of flag manifolds.

It is interesting that in every known case, {\it all}\/ solutions may be real. 
These include the classical problem of 3264 plane conics tangent to 5 plane
conics~\cite{RTV}, the 40 positions of the Stewart platform of
robotics~\cite{Dietmaier}, 
the 12 lines mutually tangent to 4 spheres~\cite{Th99}, 
the 12 rational plane cubics meeting 8 points in the plane~\cite{Kh12},
all problems of enumerating linear subspaces of a
vector space satisfying special Schubert conditions~\cite{So99}, and
certain problems of enumerating rational curves in
Grassmannians~\cite{So00}. 
These last two examples give infinitely many families of nontrivial
enumerative problems for which all solutions may be real.
They were motivated by recent, spectacular computations~\cite{FRZ,Ver98} and
a very interesting conjecture of Shapiro and Shapiro~\cite{Sottile_shapiro},
and were proved using an idea from a homotopy
continuation algorithm~\cite{HSS,HV}.

We first formalize the method of constructing real solutions introduced
in~\cite{So99,So00}, which will help extend these reality results to other
enumerative problems.
This method gives lower bounds on the maximum number of real solutions to
some enumerative problems, in the spirit of~\cite{Itenberg_Roy,St94}.
We then apply this theory to two families of enumerative problems,
one on classical (SL$_n$) flag manifolds and the other on Grassmannians of
maximal isotropic subspaces in an orthogonal vector space, showing that all
solutions may be real.
These techniques allow us to prove the opposite result---that we may have no
real solutions---for a family of enumerative problems on the Lagrangian
Grassmannian. 
Finally, we suggest a further problem to study concerning this method.

\section{Schubert Induction}

Let ${\mathbb K}$ be a field and let ${\mathbb A}^1$ be an affine 1-space over  
${\mathbb K}$.
A {\it Bruhat decomposition} of an irreducible algebraic variety $X$ defined
over ${\mathbb K}$ is a finite decomposition
$$
  X\ =\ \coprod_{w\in I} X^\circ_w
$$
satisfying the following conditions.
\begin{enumerate}
\item[(1)]
     Each stratum $X^\circ_w$ is a (Zariski) locally closed irreducible
     subvariety defined over ${\mathbb K}$ whose closure
     $\overline{X^\circ_w}$ is a union of some strata $X^\circ_v$. 
\item[(2)]
     There is a unique 0-dimensional stratum $X^\circ_{\hat{0}}$.
\item[(3)]
     For any $w,v\in I$, the intersection 
     $\overline{X^\circ_w}\cap\overline{X^\circ_v}$ is a union of some 
     strata $X^\circ_u$.
\end{enumerate}

Since $X$ is irreducible, there is a unique largest stratum
$X^\circ_{\hat{1}}$.
Such spaces $X$ include flag manifolds, where the 
$X^\circ_w$ are the Schubert cells in the Bruhat decomposition defined with
respect to a fixed flag as well as the quantum
Grassmannian~\cite{RRW98,So00,SS_SAGBI}.  
These are the only examples to which the theory developed here presently 
applies, but we expect it (or a variant) will apply to other
varieties that have such a Bruhat decomposition, particularly some spherical
varieties~\cite{Kn91} and analogs of the quantum Grassmannian for other
flag manifolds.
The key to applying this theory is to find certain geometrically
interesting families ${\mathcal Y}\to {\mathbb A}^1$ of subvarieties having
special properties with respect to the Bruhat decomposition (which we
describe below).

Suppose $X$ has a Bruhat decomposition.
Define the {\it Schubert variety} $X_w$ to be the closure of the stratum
$X^\circ_w$.
The {\it Bruhat order} on $I$ is the order induced by inclusion of Schubert
varieties: $u\leq v$ if $X_u \subset X_v$.
For flag manifolds $G/P$, these are the Schubert varieties and 
the Bruhat order on $W/W_P$; for the quantum Grassmannian, 
its quantum Schubert varieties and quantum Bruhat order.
Set $|w|:=\dim X_w$.
For flag manifolds $G/P$, if $\tau\in W$ is a minimal representative of
the coset $w\in W/W_P$ then $|w|=\ell(\tau)$, its length in the Coxeter
group $W$.

Let ${\mathcal Y}\rightarrow {\mathbb A}^1$ be a flat family of
codimension-$c$ subvarieties of $X$. 
For $s\in{\mathbb A}^1$, let $Y(s)$ be the fibre of ${\mathcal Y}$ over the
point $s$.
We say that ${\mathcal Y}$ {\it respects} the Bruhat decomposition if, for
every  $w\in I$, the (scheme-theoretic) limit 
$\lim_{s\rightarrow 0} ( Y(s)\cap X_w)$ is supported on a union of Schubert
subvarieties $X_v$ of codimension $c$ in $X_w$.
This implies that the intersection $Y(s)\cap X_w$ is proper for generic
$s\in{\mathbb A}^1$.
That is, the intersection is proper when $s$ is the generic point of the
scheme ${\mathbb A}^1$.

Given such a family, we have the cycle-theoretic equality
$$
  \lim_{s\rightarrow 0} ( Y(s)\cap X_w )\ =\ 
     \sum_{v\prec_{\mathcal Y}w} m^v_{{\mathcal Y},\,w}\, X_v\,.
$$
Here $v\prec_{\mathcal Y} w$ if $X_v$ is a component of the
support of $\lim_{s\rightarrow 0}( Y(s)\cap X_w)$, and 
the multiplicity $m^v_{{\mathcal Y},\,w}$ is the length of
the local ring of the limit scheme $\lim_{s\rightarrow 0} (Y(s)\cap X_w)$ at
the generic point of $X_v$.   
Thus, if $X$ is smooth then we have the formula
\begin{equation}\label{eq:cycle_prod}
  [X_w]\cdot[Y]\ =\ \ 
    \sum_{w\prec_{\mathcal Y}v} m^v_{{\mathcal Y},\,w}\, [X_v]\,
\end{equation}
in the Chow~\cite{Fu84a,Fu96b} or cohomology ring of $X$.
Here $[Z]$ denotes the cycle class of a subvariety $Z$, and $Y$ is any fibre
of the family ${\mathcal Y}$.
When these multiplicities $m^v_{{\mathcal Y},\,w}$ are all 1 (or 0), we call
${\mathcal Y}$ a {\it multiplicity-free family}.

A collection of families ${\mathcal Y}_1,\ldots,{\mathcal Y}_r$ 
respecting the Bruhat decomposition of $X$ is in {\it general position}
(with respect to the Bruhat decomposition) 
if, for all $w\in I$, general $s_1,\ldots,s_r\in{\mathbb A}^1$, and 
$1\leq k\leq r$, the intersection
 \begin{equation}\label{eq:intersection}
     Y_1(s_1) \cap Y_2(s_2)\cap \cdots \cap Y_k(s_k)\cap X_w
 \end{equation}
is {\it proper} in that either it is empty or else it has dimension 
$|w| -\sum_{i=1}^k c_i$, where $c_i$ is the codimension in $X$ of the fibres
of ${\mathcal Y}_i$.
Note that, more generally (and intuitively), we could require that the
intersection  
$$
  Y_{i_1}(s_{i_1}) \cap Y_{i_2}(s_{i_2})\cap \cdots \cap 
  Y_{i_k}(s_{i_k})\cap X_w
$$
be proper for any $k$-subset $\{i_1,\ldots,i_k\}$ of $\{1,\ldots,n\}$.
We do not use this added generality, although it does hold for every
application we have of this theory.
By general points $s_1,\ldots, s_k\in{\mathbb A}^1$, we mean 
general in the sense of algebraic geometry: there is a non-empty open subset
of the scheme ${\mathbb A}^k$ consisting of points $(s_1,\ldots,s_k)$ 
for which the intersection (\ref{eq:intersection}) is proper.
When $c_1+\cdots+c_k=|w|$, the
intersection~(\ref{eq:intersection}) is 0-dimensional.
Determining its degree is a problem in enumerative geometry.

We model this problem with combinatorics.
Given a collection of families ${\mathcal Y}_1,\ldots,{\mathcal Y}_r$ in
general position 
respecting the Bruhat decomposition with 
$|\hat{1}|=\dim X= c_1+\cdots+c_r$, we
construct the {\it multiplicity poset} of this enumerative
problem.
Write $\prec_i$ for $\prec_{{\mathcal Y}_i}$.
The elements of rank $k$ in the multiplicity poset are those 
$w\in I$ for which there is a chain
 \begin{equation}\label{eq:chains}
    \hat{0}\prec_1 w_1\prec_2 w_2\prec_3\cdots\prec_{k-1} w_{k-1}\prec_k
     w_k = w\,.
 \end{equation}
The cover relation between the $(i-1)$th and $i$th ranks 
is $\prec_i$.
The {\it multiplicity} of a chain~(\ref{eq:chains})
is the product of the multiplicities $m^{w_{i-1}}_{{\mathcal Y}_i,\,w_i}$ of 
the covers in that chain.
Let $\deg(w)$ be the sum of the multiplicities of all
chains~(\ref{eq:chains}) from $\hat{0}$ to $w$. 
If $X$ is smooth and $|w|=c_1+\cdots+c_k$, then $\deg(w)$ is the degree of the
intersection~(\ref{eq:intersection}), 
since it is proper, and so we have the formula~(\ref{eq:cycle_prod}).

\begin{thm}\label{thm:schubert_induction}
 Suppose $X$ has a Bruhat decomposition, 
 ${\mathcal Y}_1,\ldots,{\mathcal Y}_r$ are a collection of 
 multiplicity-free families of
 subvarieties over ${\mathbb A}^1$ in  general position,
 and each family respects this Bruhat decomposition.
 Let $c_i$ be the codimension of the fibres of ${\mathcal Y}_i$.
 \begin{enumerate}
  \item[(1)]
    For every $k$ and every $w\in I$ with $|w|=c_1+\cdots+c_k$,
    the intersection~($2$)
    is transverse for general $s_1,\ldots,s_k\in{\mathbb A}^1$ and has
    degree $\deg(w)$. 
    In particular, when ${\mathbb K}$ is algebraically closed, such an
    intersection consists of $\deg(w)$ reduced points.
  \item[(2)]
    When ${\mathbb K}={\mathbb R}$, there exist 
    real numbers $s_1,\ldots,s_r$, such that for every $k$ and every $w\in I$ with 
    $|w|=c_1+\cdots+c_k$, the intersection~($2$) is
    transverse with all points real. 
 \end{enumerate}
\end{thm}

\noindent{\bf Proof.}
For the first statement, we work in the algebraic closure of ${\mathbb K}$,
so that the degree 
of a transverse, 0-dimensional intersection is simply the number of points
in that intersection.
We argue by induction on $k$.

When $k=1$, suppose $|w|=c_1$.
Since ${\mathcal Y}_1$ is a multiplicity-free family that respects the Bruhat
decomposition, we have 
$$
  \lim_{s\rightarrow 0} ( Y_1(s)\cap X_w )\ =\ 
   m_{{\mathcal Y}_1,w}^{\hat{0}}\, X_{\hat{0}}\,,
$$
with $m_{{\mathcal Y}_1,w}^{\hat{0}}$ either 0 or 1.
Thus, for generic $s\in{\mathbb A}^1$, either $Y_1(s)\cap X_w$ is empty
or it is a single reduced point and hence transverse.
Note here that $\deg(w)= m_{{\mathcal Y}_1,w}^{\hat{0}}$.

Suppose we have proven statement (1) of the theorem for $k<l$.
Let $|w|=c_1+\cdots+c_l$.
We claim that, for generic $s_1,\ldots,s_{l-1}$, the intersection
\begin{equation}\label{eq:at_zero}
  Y_1(s_1)\cap\cdots\cap Y_{l-1}(s_{l-1})\, \cap
     \sum_{v\prec_l w} X_v
\end{equation}
is transverse and consists of $\deg(w)$ points.
Its degree is $\deg(w)$, because $\deg(w)$ satisfies the recursion
$\deg(w)=\sum_{v\prec_l w}\deg(v)$.
Transversality will follow if no two summands have a point in common.
Consider the intersection of two summands
\begin{equation}\label{eq:empty}
   Y_1(s_1)\cap\cdots\cap Y_{l-1}(s_{l-1})\cap (X_u\cap X_v)\,.
\end{equation}
Since $X_u\cap X_v$ is a union of Schubert varieties of dimensions 
less than $|w|-c_l$ and since the collection of families
${\mathcal Y}_1,\ldots,{\mathcal Y}_{l-1}$ is  in general
position, it folows that~(\ref{eq:empty}) is empty for generic
$s_1,\ldots,s_{l-1}$, which proves transversality. 
Consider now the family defined by $Y_l(s)\cap X_w$ for $s$ generic.
Since $\sum_{v\prec_l w} X_v$ is the 
fibre of this family at $s=0$ and since 
the intersection~(\ref{eq:at_zero}) is transverse and consists of
$\deg(w)$ points, for generic 
$s_l\in{\mathbb A}^1$ the intersection
 \begin{equation}\label{eq:l-intersection}
     Y_1(s_1) \cap Y_2(s_2)\cap \cdots 
      \cap Y_{l-1}(s_{l-1})\cap Y_l(s_l)\cap X_w
 \end{equation}
is transverse and consists of $\deg(w)$ points. 
\smallskip

For statement (2) of the theorem, we inductively construct 
real numbers $s_1,\ldots,s_r$ having the properties that:
(a) for any $w\in I$ and $k$ with 
$|w|=c_1+\cdots+c_k$, the intersection~(\ref{eq:intersection}) is transverse
with all points real; and (b) that if $|w|<c_1+\cdots+c_k$,
then~(\ref{eq:intersection}) is empty. 
Suppose $|w|=c_1$.
Since for general $s\in{\mathbb R}$ the intersection
$X_w\cap Y_1(s)$ is either empty or consists of a single reduced 
point, we may select a general $s\in{\mathbb R}$ with the additional
property that if $|v|<c_1$ then $Y_1(s)\cap X_v$ is empty.

Suppose now that we have constructed $s_1,\ldots,s_{l-1}\in{\mathbb R}$ such
that (a) if $|v|=c_1+\cdots+c_{l-1}$ then the intersection
$Y_1(s_1)\cap \cdots\cap Y_{l-1}(s_{l-1})\cap X_v$ is transverse 
with all points real, and (b) if $|v|<c_1+\cdots+c_{l-1}$, then
this intersection is empty.
Let $|w|=c_1+\cdots+c_l$.
Then the intersection~(\ref{eq:at_zero}) is transverse with all points real.
Thus there exists $\epsilon_w>0$ such that if
$0<s_l\leq\epsilon_w$, then the intersection~(\ref{eq:l-intersection}) is
transverse with all points real.
Set $s_l=\min\{\epsilon_w: |w|=c_1+\cdots+c_l\}$.
Since it is an open condition (in the usual topology) on the $l$-tuple 
$(s_1,\ldots,s_l)\in{\mathbb R}^l$ for the
intersection~(\ref{eq:l-intersection}) to be transverse with all points real
and since there are finitely many $w\in I$, we may (if necessary) choose
a nearby $l$-tuple of points  
such that, if $|w|<c_1+\cdots+c_l$, then the 
intersection~(\ref{eq:l-intersection}) is empty. 
\qed\smallskip

\begin{rem}
The statement and proof of Theorem~\ref{thm:schubert_induction} 
are a generalization of the main results of~\cite[Thm.~1]{So99}
and~\cite[Thms.~3.1 and~3.2]{So00} and they constitute a stronger version
of the theory presented in~\cite{So97b}.
(Part 1 generalizes~\cite[Thm.~8.3]{EH83}). 
We call this method of proof {\it Schubert induction}.
The proof of the second statement is based upon the fact that small (real)
perturbations of a transverse intersection preserve transversality as well as
the number of real and complex points in that intersection.
In principle, this leads to an optimal  numerical homotopy continuation
algorithm for finding all complex points in the
intersection~(\ref{eq:intersection}). 
A construction and correctness proof of such an algorithm could be modeled
on the Pieri homotopy algorithm of~\cite{HSS,HV}.
\end{rem}

\begin{rem}
The first statement of Theorem~\ref{thm:schubert_induction} gives an
elementary proof of generic transversality for some enumerative problems
involving multiplicity-free families.
In characteristic 0, it is an alternative to Kleiman's Transversality
Theorem~\cite{MR50:13063} and could provide a basis to prove generic 
transversality in arbitrary characteristic,
extending the result in~\cite{So97a} that the intersection of general
Schubert varieties in a Grassmannian of 2-planes is generically transverse
in any characteristic.
It also provides a proof that $\deg(w)$ is the intersection number---without
using Chow or cohomology rings, the traditional tool in enumerative 
geometry. 
\end{rem}

\begin{rem}
If the families ${\mathcal Y}_i$ are not multiplicity-free, then we can
prove a lower bound on the maximum number of real solutions.
A (saturated) chain~(\ref{eq:chains}) in the multiplicity poset is
{\it odd}\/ if it has odd multiplicity.
Let $\mbox{\rm odd}(w)$ count the odd chains from $\hat{0}$ to $w$ in the
multiplicity poset.

\begin{thm}\label{thm:bound}
 Suppose $X$ has a Bruhat decomposition,
 ${\mathcal Y}_1,\ldots,{\mathcal Y}_r$ are a collection of families of
 subvarieties over ${\mathbb A}^1$ in  general position,
 and each family respects this Bruhat decomposition.
 Let $c_i$ be the codimension of the fibres of ${\mathcal Y}_i$.
 \begin{enumerate}
  \item[(1)]
    Suppose ${\mathbb K}$ is algebraically closed.
    For every $k$, every $w\in I$ with $|w|=c_1+\cdots+c_k$, and general
     $s_1,\ldots,s_k\in{\mathbb A}^1$, 
    the $0$-dimensional intersection~($2$) has degree
    $\deg(w)$.
  \item[(2)]
    When ${\mathbb K}={\mathbb R}$, there exist real numbers
    $s_1,\ldots,s_r$ such that for every $k$, every $w\in I$ with 
    $|w|=c_1+\cdots+c_k$, the intersection~($2$) is
    $0$-dimensional and has at least $\mbox{\rm odd}(w)$ real points. 
 \end{enumerate}
\end{thm}

\noindent{\bf Sketch of Proof. }
For the first statement, the same arguments as in the proof of
Theorem~\ref{thm:schubert_induction}  suffice if we replace the phrase
``transverse and consists of $\deg(w)$ points'' throughout by 
``proper and has degree $\deg(w)$''.
For statement (2) of the theorem, observe that a point in the
intersection  $Y_1(s_1)\cap \cdots\cap Y_{l-1}(s_{l-1})\cap X_v$ 
becomes $m^v_{{\mathcal Y}_l,\,w}$ points counted with multiplicity
in~(\ref{eq:l-intersection}), when $s_l$ is a small real number. 
If this multiplicity $m^v_{{\mathcal Y}_l,\,w}$ is odd and the original
point was real, then at least one of these $m^v_{{\mathcal Y}_l,\,w}$ 
points are real.
\qed\medskip

The lower bound of Theorem~\ref{thm:bound} is the analog of the bound for
sparse polynomial systems in terms of alternating mixed
cells~\cite{Itenberg_Roy,PS94,St98}.   
Like that bound, it is not sharp~\cite{LW98,St98}.  
We give an example using the notation of Section~2.
The Grassmannian of 3-planes in ${\mathbb C}^7$
has a Bruhat decomposition indexed by triples 
$1\leq \alpha_1<\alpha_2<\alpha_3\leq 7$ of
integers. 
Let $r=4$ and suppose that each family ${\mathcal Y}_i$ is the family 
of Schubert varieties $X_{357}\Fdot(s)$,
where $\Fdot(s)$ is the flag of subspaces osculating a real rational normal
curve.
In~\cite[Thm.~3.9(iii)]{Sottile_shapiro} it is proven that if $s,t,u,v$
are distinct real points, then 
$$
  Y(s)\cap Y(t)\cap Y(u)\cap Y(v)
$$
is transverse and consists of eight real points.
However, there are  fice chains in the  multiplicity poset; four of them odd and
one of multiplicity 4.
In Figure~\ref{fig:chains}, we show the Hasse diagram of this multiplicity
poset, indicating multiplicities greater than 1.
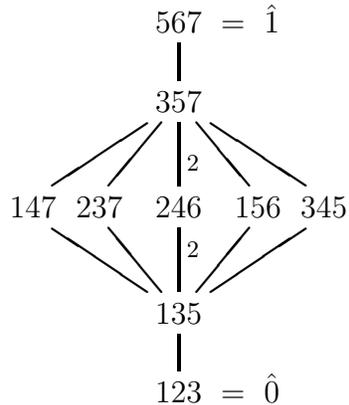
\begin{figure}[htb]
$$
\begin{picture}(130, 150)(5,0)
\thicklines
\put(70,122){\line(0,1){14}}
\put(70,106){\line(0,-1){24}}\put(73,88){\scriptsize $2$}
\put(58,106){\line(-3,-2){36}}\put(63.3,106){\line(-5,-6){20}}
\put(82,106){\line( 3,-2){36}}\put(76.6,106){\line( 5,-6){20}}
\put(58, 42){\line(-3, 2){36}}\put(63.3, 42){\line(-5, 6){20}}
\put(82, 42){\line( 3, 2){36}}\put(76.6, 42){\line( 5, 6){20}}
\put(70, 42){\line(0,1){24}}\put(73,55){\scriptsize $2$}
\put(70, 12){\line(0,1){14}}

\put( 61,140){$567\ =\ \hat{1}$}
\put( 61,110){$357$}
\put(  6, 70){$147$}\put( 31, 70){$237$}\put( 61, 70){$246$}
\put( 91, 70){$156$}\put(116, 70){$345$}
\put( 61, 30){$135$}\put( 61,  0){$123\ =\ \hat{0}$}
\end{picture}
$$
 \caption{The multiplicity poset\label{fig:chains}} 
\end{figure}

Despite this lack of sharpness, Theorem~\ref{thm:bound} gives new results
for the Grassmannian.
In~\cite{EH87}, Eisenbud and Harris show that families of Schubert
subvarieties of a Grassmannian defined by flags of subspaces osculating a
rational normal curve respect the Bruhat decomposition given by any such
osculating flag, and any collection is in general position.
Consequently, given a collection of these families with 
$\mbox{\rm odd}(w)>0$, it follows that $\mbox{\rm odd}(w)$
is a nontrivial lower bound (new if the
Schubert varieties are not special Schubert varieties) on the
number of real points in such a 0-dimensional intersection of these Schubert
varieties.

For example, in the Grassmannian of 3-planes in ${\mathbb C}^{r+3}$, let
$Y(s)$ be the Schubert variety consisting of 3-planes having nontrivial
intersection with $F_{r-1}(s)$ and whose linear span with $F_{r+1}(s)$
is not all of ${\mathbb C}^{r+3}$.
(Here, $F_i(s)$ is the $i$-dimensional subspace osculating a real rational
normal curve $\gamma$ at the point $\gamma(s)$.)
This Schubert variety has codimension 3.
Consider the enumerative problem given by intersecting $r$ of these
Schubert varieties.
Table~\ref{table:one} gives both the number of solutions ($\deg(\hat{1})$)
\begin{table}[htb]
 \begin{tabular}{|c||c|c|c|c|c|c|c|c|c|c|c|c|}\hline\hline
  $r$& 2& 3& 4 & 5 & 6& 7 & 8 &9 &10 &11\\\hline
  $\deg(\hat{1})$& \rule{0pt}{13pt}
       1&2&8&32&145&702&3598&19,280&107,160&614,000\\\hline
  $\mbox{\rm odd}(\hat{1})$&\rule{0pt}{13pt}
       1&0&4& 6& 37&116& 534& 2128&  9512& 41,656\\\hline\hline
 \end{tabular}\bigskip
 \caption{Number of solutions and odd chains\label{table:one}}
\end{table}
and the number of odd chains ($\mbox{\rm odd}(\hat{1})$) in the multiplicity
poset for $r=2,3,\ldots,11$. 
The case $r=4$ we have already described.
The conjecture of Shapiro and Shapiro~\cite{Sottile_shapiro}
asserts that all solutions for any $r$-tuple of distinct real points will be
real, which is stronger than the consequence of
Theorem~\ref{thm:bound}  that there is some $r$-tuple of real points for
which there will be at least as many  real solutions as odd chains. 
\end{rem}

\begin{rem}
The requirement that there be a unique 0-dimensional stratum in a
Bruhat decomposition may be relaxed.
We could allow several 0-dimensional strata $X_z$ for $z\in Z$,
each consisting of a single ${\mathbb K}$-rational point.
This is the case for toric varieties~\cite{Fulton_toric} and more generally
for spherical varieties~\cite{Kn91}.

If we define the multiplicity poset as before, then $Z$ indexes 
its minimal elements.
We define the intersection number $\deg(w)$ and the bound 
$\mbox{\rm odd}(w)$ using chains
$$
  z \prec_1 w_1 \prec_2 w_2 \prec_3 \cdots \prec_x w_k \ =\ w\,
    \quad  \mbox{ with }\quad z\in Z\,.
$$
Then almost the same proof as we gave for
Theorem~\ref{thm:schubert_induction} proves the same statement in this new
context.  
We do not yet know of any applications of this extension of
Theorem~\ref{thm:schubert_induction}, but we expect that some will be found.
\end{rem}

\section{The Classical Flag Manifolds}
Fix integers $n\geq m>0$ and a sequence
${\bf d} :0<d_1<\cdots<d_m<n$ of integers. 
A {\it partial flag of type ${\bf d}$} is a sequence of linear subspaces
$$
  E_{d_1}\ \subset\ E_{d_2}\ \subset\ \cdots\ 
     \subset\ E_{d_m}\ \subset\ {\mathbb C}^n
$$
with $\dim E_i=d_i$ for each $i=1,\ldots,m$.
The {\it flag manifold}\/ $\Fla$ is the collection of all
partial flags of type ${\bf d}$. 
This manifold is the homogeneous space
SL$(n,{\mathbb C})/P_{\bf d}$, where $P_{\bf d}$ is the parabolic subgroup of
SL$(n,{\mathbb C})$ defined by the simple roots {\it not} indexed by
$\{d_1,\ldots,d_m\}$.
See~\cite{Bo91} or~\cite{Fu97} for further
information on partial flag varieties.

A fixed {\it complete flag} $\Fdot$ 
($F_1\subset\cdots\subset F_n={\mathbb C}^n$ with $\dim F_i=i$)
induces a Bruhat decomposition of $\Fla$ 
\begin{equation}\label{eq:Fla_schubert}
  \Fla\ =\ \coprod X^\circ_w\Fdot
\end{equation}
indexed by those permutations $w=w_1\ldots w_n$ in the symmetric group
${\mathcal S}_n$ whose descent set $\{i\mid w_i>w_{i+1}\}$ is a
subset of $\{d_1,\ldots,d_m\}$.
Write $I_{{\bf d}}$ for this set of permutations.
Then $|w|=\ell(w)$, as $I_{{\bf d}}$ is the set of minimal coset
representatives for $W_{P_{\bf d}}$.
The Schubert variety $X_w\Fdot$ is the closure of the Schubert cell 
$X^\circ_w\Fdot$.

Fix any real rational normal curve $\gamma\colon{\mathbb C}\to{\mathbb C}^n$,
which is a map given by $\gamma\colon s\mapsto(p_1(s),\ldots,p_n(s))$, where
$p_1,\ldots,p_n$ 
are a basis for the space of real polynomials of degree less than $n$.
All real rational normal curves are isomorphic by a real linear
transformation.
For any $s\in {\mathbb C}$, let $\Fdot(s)$ be the complete flag of
subspaces osculating the curve $\gamma$ at the point
$\gamma(s)$.
The dimension-$i$ subspace $F_i(s)$ of $\Fdot(s)$ is 
the linear span of the vectors 
$\gamma(s)$ and $\gamma'(s):=\frac{d}{ds}\gamma(s),\ldots,\gamma^{(i-1)}(s)$.

For each $i=1,\ldots,m$, we have {\it simple Schubert variety}
$X_i\Fdot$ of $\Fla$.
Geometrically, 
$$
   X_i\Fdot\ :=\ 
   \{\Edot\in \Fla\mid E_{d_i}\cap F_{n-d_i}\neq\{0\}\}\,.
$$
We call these ``simple'' Schubert varieties, for they give simple 
(codimension-1) conditions on partial flags in $\Fla$.
Let ${\mathcal X}_i\to{\mathbb A}^1$ be the family 
whose fibre over $s\in{\mathbb A}^1$ is $X_i\Fdot(s)$. 
We study these families.

\begin{thm}\label{thm:family}
Let ${\bf d}=0<d_1<\cdots<d_m<n$ be a sequence of integers.
For any $i=1,\ldots, m$, the family 
${\mathcal X}_i\to{\mathbb A}^1$ of simple Schubert varieties
is a multiplicity-free family that respects the Bruhat decomposition of\/
$\Fla$ given by the flag $\Fdot(0)$.

Any collection of these families of simple Schubert varieties is in general
position. 
\end{thm}

We shall prove Theorem~\ref{thm:family} shortly.
First, by Theorem~\ref{thm:schubert_induction}, we deduce the following 
corollary.

\begin{cor}\label{cor:flag-reality}
Let $w\in I_{{\bf d}}$ and set $r:=|w|=\dim X_w$.
Then, for any list of numbers $i_1,\ldots,i_r\in\{1,\ldots,m\}$,
there exist real numbers $s_1,\ldots,s_r$ such that 
\begin{equation}\label{eq:real-trans}
  X_w\Fdot(0)\cap X_{i_1}\Fdot(s_1)\cap\cdots\cap
   X_{i_r}\Fdot(s_r)
\end{equation}
is transverse and consists only of real points.
\end{cor}

This corollary generalizes the intersection of the main results of~\cite{So99}
and~\cite{So00}, which is the case of Corollary~\ref{cor:flag-reality} for
Grassmannians (${\bf d}=d_1$ has only a single part).
This result also extends (part of) Theorem~13 in~\cite{So97b}, which states
that, if ${\bf d}=2<n-2$ and $i_1,\ldots,i_r$ are any numbers from
$\{2,n-2\}$  ($r=\dim \Fla=4n-12$), then
there exist real flags $F_\bullet^1,\ldots,F_\bullet^r$ such that 
$$
  X_{i_1}F_\bullet^1\cap\cdots\cap  X_{i_r}F_\bullet^r
$$
is transverse and consists only of real points.
\medskip

We recall some additional facts about the
cohomology of the partial flag manifolds $\Fla$.
Each stratum $X^\circ_w\Fdot$ is isomorphic to ${\mathbb C}^{|w|}$ and 
the Bruhat decomposition~(\ref{eq:Fla_schubert}) is a cellular decomposition
of $\Fla$ into even- (real) dimensional cells.
Let $\sigma_w$ be the cohomology class Poincar\'e dual to the fundamental
(homology) cycle of the Schubert variety $X_w\Fdot$.
Then these Schubert classes $\sigma_w$ provide a basis for the integral
cohomology ring $H^*(\Fla,{\mathbb Z})$ with 
$\sigma_w\in  H^{2c(w)}(\Fla,{\mathbb Z})$, where $c(w)$ is the 
complex codimension of $X_w\Fdot$ in $\Fla$.

Let $\tau_i$ be the class of the simple Schubert variety $X_i\Fdot$.
There is a simple formula due to
Monk~\cite{Monk} and Chevalley~\cite{Ch91} expressing the product
$\sigma_w\cdot\tau_i$ in terms of the basis of Schubert classes.
Let $w\in I_{{\bf d}}$.
Then 
$$
  \sigma_w\cdot \tau_i \ =\ 
  \sum \sigma_{w(j,k)}\,,
$$
where $(j,k)$ is a transposition; the sum is over all 
$j\leq d_i<k$, where
 \begin{enumerate}
  \item[(1)] $w_j>w_k$ and
  \item[(2)] if $j<l<k$ then either $w_l>w_j$ or else $w_k>w_l$.
 \end{enumerate}
Write $w(j,k)\lessdot_i w$ for such $w(j,k)$.
Note that, if $w\in I_{{\bf d}}$, then so is any 
$v\in{\mathcal S}_n$ with $v\lessdot_i w$ for any $i=1,\ldots,m$.

Let $Gr(d_i)$ be the Grassmannian of $d_i$-dimensional subspaces of
${\mathbb C}^n$.
The association $\Edot\mapsto E_{d_i}$ induces a projection 
$\pi_i:\Fla\to Gr(d_i)$.
The Grassmannian has a Bruhat decomposition 
$$
  Gr(d_i)\ =\ \coprod \Omega^\circ_\alpha\Fdot
$$
indexed by increasing sequences $\alpha$ of length $d_i$, 
$1\leq\alpha_1<\alpha_2<\cdots<\alpha_{d_i}\leq n$, with the Bruhat order
given by componentwise comparison.
Such an increasing sequence can be uniquely completed to a permutation
$w(\alpha)$ whose only descent is at $d_i$.
The map $\pi_i$ respects the two Bruhat decompositions in that
$\pi_i^{-1}(\Omega_\alpha)=X_{w(\alpha)}\Fdot$ and 
$\pi_i(X_w\Fdot) = \Omega_{\alpha(w)}\Fdot$, 
where $\alpha(w)$ is the sequence obtained by writing
$w_1,\ldots,w_{d_i}$ in increasing order.
Thus, if $\beta<\alpha(w)$, then $X_w\Fdot\cap\pi_i^{-1}\Omega_\beta\Fdot$ is
a union of proper Schubert subvarieties of $X_w\Fdot$.\smallskip

The Grassmannian has a distinguished simple Schubert variety
$$
  \Upsilon\Fdot\ =\ 
   \{ E\in Gr(d_i)\mid E\cap F_{n-d_i} \neq\{0\}\}\,.
$$
This shows 
$X_i\Fdot\ =\ \pi_i^{-1}(\Upsilon\Fdot)$.
We have $\Upsilon\Fdot=\Omega_{(n-d_i,n-d_i+2,\ldots,n)}\Fdot$.

We need the following useful fact about the families 
${\mathcal X}_w\to{\mathbb A}^1$.

\begin{lemma}\label{lem:infint}
For any $w\in I_{{\bf d}}$, we have 
${\displaystyle \bigcap_{s\in{\mathbb A}^1}X_w\Fdot(s)\ =\ \emptyset}$.
\end{lemma}

\noindent{\bf Proof.} 
Any Schubert variety $X_w\Fdot$ is a subset of some simple Schubert variety 
$X_i\Fdot=\pi_i^{-1}\Upsilon\Fdot$.
Thus it suffices to prove the lemma for the simple Schubert
varieties $\Upsilon\Fdot(s)$ of a Grassmannian.
But this is simply a consequence of~\cite[Thm.~2.3]{EH83}.
\qed\medskip

\noindent{\bf Proof of Theorem~\ref{thm:family}.}
For any $w\in I_{{\bf d}}$, we consider the scheme-theoretic limit
$\lim_{s\to 0}(X_w\Fdot(0)\cap X_i\Fdot(s))$.
Since $X_i\Fdot\ =\ \pi_i^{-1}(\Upsilon\Fdot)$, for any 
$s\in{\mathbb C}$ we have
$$
  X_w\Fdot(0)\cap X_i\Fdot(s)\ =\ 
     X_w\Fdot(0)\cap \pi_i^{-1}\left(
       \Omega_{\alpha(w)}\Fdot(0)\cap \Upsilon\Fdot(s)\right)\,,
$$
since $\pi_i X_w\Fdot(0)=\Omega_{\alpha(w)}\Fdot(0)$.
Thus, set-theoretically we have 
$$
  \lim_{s\to 0}\left(X_w\Fdot(0)\cap X_i\Fdot(s)\right)\ \subset\ 
     X_w\Fdot(0)\cap \pi_i^{-1}\left( \lim_{s\to 0}\left(
       \Omega_{\alpha(w)}\Fdot(0)\cap \Upsilon\Fdot(s)\right)\right)\,.
$$
But this second limit is $\bigcup_{\beta<\alpha(w)} \Omega_\beta\Fdot$
by~\cite[Thm.~8.3]{EH83}. 
Thus
\begin{eqnarray*}
  \lim_{s\to 0}\left(X_w\Fdot(0)\cap X_i\Fdot(s)\right) &\subset&
  X_w\Fdot(0)\cap \pi_i^{-1}\left(
         \bigcup_{\beta<\alpha(w)} \Omega_\beta\Fdot(0)\right)\\
  &\subset& \bigcup_{v\lessdot w} X_v\Fdot(0)\,,
\end{eqnarray*}
set-theoretically.

Since the limit scheme 
$\lim_{s\to 0}\left(X_w\Fdot(0)\cap X_i\Fdot(s)\right)$ is supported on 
this union of proper Schubert subvarieties of $X_w\Fdot(0)$ and has
dimension at least $\dim X_w\Fdot(0)-1$, its support must be a union of
codimension-1 Schubert subvarieties of $X_w\Fdot(0)$. 
Hence the family ${\mathcal X}_i\to{\mathbb A}^1$ respects the Bruhat
decomposition, and we have
$$
  \lim_{s\to 0}\left(X_w\Fdot(0)\cap X_i\Fdot(s)\right)\ =\ 
  \sum_{v\lessdot w} m^v_{i,\,w}\, X_v\Fdot(0)\,.
$$
thus $\sigma_w\cdot\tau_i=\sum_{v\lessdot w} m^v_{i,\,w}\sigma_v$
in the Chow ring.
Since the Schubert classes $\sigma_v$ are linearly independent in the Chow
ring, these multiplicities are either 0 or 1 by Monk's formula, and they are 1
precisely when $v\lessdot_i w$.
Thus the family ${\mathcal X}_i\to{\mathbb A}^1$ is multiplicity-free,
and we have proven the first statement of Theorem~\ref{thm:family}.
\smallskip

To complete the proof, 
let ${\mathcal X}_{i_1},\ldots,{\mathcal X}_{i_r}$ be a collection of
families of simple 
Schubert varieties defined by the flags $\Fdot(s)$.
We show that this collection is in general
position with respect to the Bruhat decomposition defined by the flag
$\Fdot(0)$. 
If not, then there is some index $w$ and integer $k$ with $k$ minimal such
that, for general $s_1,\ldots, s_k\in{\mathbb C}$,
\begin{equation}\label{eq:genint}
  X_w\Fdot(0)\cap X_{i_1}\Fdot(s_1)\cap\cdots\cap X_{i_{k-1}}\Fdot(s_{k-1})
\end{equation}
has dimension $|w|-k+1$, but 
$$
  X_w\Fdot(0)\cap X_{i_1}\Fdot(s_1)\cap\cdots\cap X_{i_k}\Fdot(s_k)
$$ 
has dimension exceeding $|w|-k$.
Hence its dimension is $|w|-k+1$.
But then, for general $s\in{\mathbb C}$, some component
of~(\ref{eq:genint}) lies in $X_{i_k}\Fdot(s)$, which implies that this
component lies in $X_{i_k}\Fdot(s)$ for all $s\in{\mathbb C}$, 
contradicting Lemma~\ref{lem:infint}.
\qed\medskip

The previous paragraph provides a proof of the following useful lemma.

\begin{lemma}\label{lem:simple-genpos}
Suppose a variety $X$ has a Bruhat decomposition. 
Let ${\mathcal Y}_1,\ldots,{\mathcal Y}_r$ be a collection of codimension-$1$
families in $X$, each of which respects this Bruhat decomposition.
If each family ${\mathcal Y}_i\to{\mathbb A}^1$ satisfies
$$
  \bigcap_{s\in{\mathbb A}^1} Y_i(s)\ =\ \emptyset\,,
$$
then the collection of families ${\mathcal Y}_1,\ldots,{\mathcal Y}_r$
is in general
position.
\end{lemma}

A fruitful question is to ask how much freedom we have to select the
real numbers $s_1,\ldots,s_r$ of Corollary~\ref{cor:flag-reality} so that
all the points of the intersection~(\ref{eq:real-trans}) are real.
In 1995, Boris Shapiro and Michael Shapiro conjectured that we have almost
complete freedom:
For generic real numbers $s_1,\ldots,s_r$, all points
of~(\ref{eq:real-trans}) are real. 
This remarkable conjecture is false in a very interesting way.

\begin{ex}\label{ex:NC}
Let $n=5$ and ${\bf d}:2<3$ so that $\Fla$ is the manifold of flags
$E_2\subset E_3\subset{\mathbb C}^5$.
This 8-dimensional flag manifold has two types of simple Schubert varieties
$X_i\Fdot$ for $i=2,3$, where $X_i\Fdot$ consists of those flags 
$E_2\subset E_3$ with $E_i\cap F_{5-i}\neq\{0\}$.
Write $X_i(s)$ for $X_i\Fdot(s)$.
A calculation (using Maple and Singular~\cite{SINGULAR}) 
shows that 
$$
  X_2(-8)\cap   X_3(-4)\cap   X_2(-2)\cap   X_3(-1)\cap 
  X_2(1)\cap   X_3(2)\cap   X_2(4)\cap   X_3(8)
$$
is transverse and consists of twelve points, none of which are real.
\end{ex}

Despite this counterexample, quite a lot may be salvaged from the conjecture
of Shapiro and Shapiro.
When the partial flag manifold $\Fla$ is a Grassmannian, there are no known
counterexamples, many enumerative problems, and choices of 
real numbers $s_1,\ldots,s_r$ for which all solutions are 
real~\cite{Sottile_shapiro}; in~\cite{EG00}, the conjecture is proven 
for any Grassmannian of 2-planes. 
The general situation seems much subtler.
In our counterexample, the points $\{-8,-2,1,4\}$  at which we
evaluate $X_2$ alternate with the points $\{-4,-1,2,8\}$ at which we
evaluate $X_3$.
If, however, we evaluate $X_2$ at points $s_1,\ldots,s_4$ and $X_3$ at points
$s_5,\ldots,s_8$ with $s_1<s_2<\cdots<s_8$, then we know of no
examples with any points of intersection not real.
We have checked this for all 24,310 subsets of eight numbers from 
$$
  \{-6,-5,-4,-3,-2,-1,1,2,3,5,7,11,13,17,19,23,29\}\,.
$$
On the other hand, if we evaluate $X_2$ at any four of the eight numbers
$$
  \{1,2, 3^2, 4^3, 5^4, 6^5, 7^6, 8^7\}
$$ 
and $X_3$ at the other four numbers,
then all twelve points of intersection are real.

\section{The Orthogonal Grassmannian}
Let $V$ be a vector space equipped with a nondegenerate symmetric bilinear
form $\Span{\cdot,\cdot}$.
A subspace $H\subset V$ is {\it isotropic} if the restriction of the form
to $H$ is identically zero.
Isotropic subspaces have dimension at most half that of $V$.
The {\it orthogonal Grassmannian} is the collection of all isotropic
subspaces of $V$ with maximal dimension.
If the dimension of $V$ is even, then the orthogonal Grassmannian has two
connected components, and each is isomorphic to the orthogonal
Grassmannian for a generic hyperplane section of $V$;
the isomorphism is given 
by intersecting with that hyperplane.
Thus, it suffices to consider only the case when the dimension of $V$ is odd.

When $V$ has dimension $2n+1$, a maximal isotropic subspace $H$ of $V$ has
dimension $n$, and we write $OG(n)$ for this orthogonal Grassmannian.
To ensure that $OG(n)$ has ${\mathbb K}$-rational points, we assume that 
$V$ has a ${\mathbb K}$-basis $e_1,\ldots,e_{2n+1}$, for which our form is
\begin{equation}\label{eq:orth-form}
  \left\langle \sum x_ie_i,\ \sum y_j e_j\right\rangle\ =\ 
   \sum x_i y_{2n+2-i}\,.
\end{equation}
Then $OG(n)$ is a homogeneous space of the (split) special orthogonal group
$SO(2n+1,{\mathbb K})=\mbox{Aut}(V,\Span{\cdot,\cdot})$.
This algebraic manifold has dimension $\binom{n+1}{2}$.

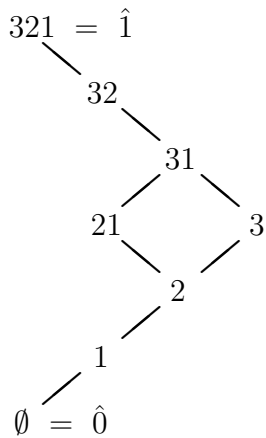
\begin{figure}[htb]
$$
\begin{picture}(95,163)(-1,-3)
\thicklines
                            \put(25,136){\line(-6,5){14}}
                            \put(55,111){\line(-6,5){14}}
\put(41, 86){\line(6,5){14}}\put(85,86){\line(-6,5){14}}
\put(71, 61){\line(6,5){14}}\put(55,61){\line(-6,5){14}}
\put(41, 36){\line(6,5){14}}
\put(11, 11){\line(6,5){14}}

\put(-2,150){$321\ =\ \hat{1}$}
\put(27.5,125){$32$}
\put(57,100){$31$}
\put(29, 75){$21$}\put(89,75){$3$}
\put(59, 50){$2$}
\put(30, 25){$1$}
\put( 0,  0){$\emptyset\ =\ \hat{0}$}
\end{picture}
$$
 \caption{The Bruhat order for $OG(3)$.\label{fig:BO}} 
\end{figure}

An {\it isotropic flag} is a complete flag $\Fdot$ of $V$ such that 
(a) $F_n$ is isotropic and (b) for every
$i> n$, $F_i$ is the annihilator of $F_{2n+1-i}$,
that is, $\Span{F_{2n+1-i},\,F_i}\equiv 0$.
An isotropic flag induces a Bruhat decomposition 
$$
  OG(n)\ =\ \coprod X_\lambda^\circ \Fdot
$$
indexed by decreasing
sequences $\lambda$ of positive integers
$n\geq\lambda_1>\cdots>\lambda_l>0$, called {\it strict partitions}. 
Let $SP(n)$ denote this set of strict partitions.
The Schubert variety $X_\lambda\Fdot$ is the closure of 
$X_\lambda^\circ\Fdot$ and has dimension 
$|\lambda|:=\lambda_1+\cdots+\lambda_l$.
The Bruhat order is given by componentwise comparison:
$\lambda\geq\mu$ if $\lambda_i\geq\mu_i$ for all $i$ with both
$\lambda_i,\mu_i>0$. 
Figure~\ref{fig:BO} illustrates this Bruhat order when $n=3$.

The unique simple Schubert variety of $OG(n)$ is (set-theoretically)
$$
  Y\!\Fdot\ :=\ \{H\in OG(n)\mid H\cap F_{n+1}\neq \{0\}\}\,.
$$
Thus $Y\!\Fdot$ is the set-theoretic intersection of $OG(n)$ with the simple
Schubert variety $\Upsilon\Fdot$  of the ordinary Grassmannian $Gr(n)$ of
$n$-dimensional subspaces of $V$.
The multiplicity of this intersection is 2 (see~\cite[p.~68]{FuPr}).
We have $Y\!\Fdot=X_{(n,n-1,\ldots,2)}\Fdot$.
The Bruhat orders of these two Grassmannians ($OG(n)$ and $Gr(n)$) are
related. 

\begin{lemma}\label{lem:BO-rel}
Let $\Fdot$ be a fixed isotropic flag in $V$.
Then every Schubert cell $X^\circ_\lambda\Fdot$ of $OG(n)$ lies in a unique
Schubert cell $\Omega^\circ_{\alpha(\lambda)}\Fdot$ of $Gr(n)$.
Moreover, for any strict partition $\lambda$, we have the set-theoretic
equality
$$
  X_\lambda\Fdot\cap 
  \bigcup_{\beta\lessdot\alpha(\lambda)}\Omega_\beta\Fdot\ =\ 
  \bigcup_{\mu\lessdot\lambda} X_\mu\Fdot\,.
$$
\end{lemma}

Let $\tau$ be the cohomology class dual to the fundamental cycle of
$Y\!\Fdot$, and let $\sigma_\lambda$ be the class dual to the fundamental
cycle of $X_\lambda\Fdot$. 
The Chevalley formula for $OG(n)$ is
$$
  \sigma_\lambda\cdot \tau\ = \ \sum_{\mu\lessdot\lambda}\sigma_\mu\,,
$$
which is free of multiplicities.

Let ${\mathbb K}={\mathbb C}$.
As in Section 2, we study families of Schubert varieties defined by flags
$\Fdot(s)$ of isotropic subspaces osculating a real rational normal curve
$\gamma:{\mathbb C}\to V$ at $\gamma(s)$.
With our given form $\Span{\cdot,\cdot}$ and basis $e_1,\ldots,e_{2n+1}$,
one choice for a real rational normal
curve $\gamma$ whose flags of osculating subspaces are isotropic is 
$$
  \gamma(s)\ =\ \left(
      1,\,s,\,\frac{s^2}{2},\,\ldots,\, \frac{s^n}{n!},\,
     -\frac{s^{n+1}}{(n+1)!},\,\frac{s^{n+2}}{(n+2)!},\,
      \ldots,\,(-1)^n\frac{s^{2n}}{(2n)!}\right)\,.
$$

\begin{thm}\label{thm:orth-families}
The family ${\mathcal Y}\to{\mathbb A}^1$ of simple Schubert
varieties $Y\!\Fdot(s)$ is multiplicity-free and respects the Bruhat
decomposition of $OG(n)$ induced by the flag $\Fdot(0)$.

Any collection of these families of simple Schubert varieties is 
in general position.
\end{thm}

We omit the proof of this theorem, which is nearly identical to the proof of
Theorem~\ref{thm:family}.
By Theorem~\ref{thm:schubert_induction}, we deduce the following 
corollary.

\begin{cor}\label{cor:orth-reality}
Let $\lambda\in SP(n)$.
Then there exist real numbers $s_1,\ldots,s_{|\lambda|}$ such that 
\begin{equation}\label{eq:orth-int}
  X_\lambda\Fdot(0)\cap Y\!\Fdot(s_1)\cap\cdots\cap
   Y\!\Fdot(s_{|\lambda|})
\end{equation}
is transverse and consists only of real points.
\end{cor}

By Theorem~\ref{thm:orth-families} and the Chevalley formula, for a strict
partition $\lambda$ and 
general complex numbers $s_1,\ldots,s_{|\lambda|}$, 
the intersection~(\ref{eq:orth-int}) is transverse and  
consists of $\deg(\lambda)$ points, where $\deg(\lambda)$ is the
number of chains in the Bruhat order from $0=\hat{0}$ to $\lambda$.

As in Section 2, we may ask how much freedom we have to select the
real numbers $s_1,\ldots,s_{|\lambda|}$ of
Corollary~\ref{cor:orth-reality} 
so that all the points of the intersection~(\ref{eq:orth-int}) are real.
When $n=3$ and $\lambda=\hat{1}$ (Figure~\ref{fig:BO} shows that
$|\hat{1}|=6$ and $\deg(\hat{1})=2$), the discriminant of
a polynomial formulation of this problem is
$$
  \sum_{w\in {\mathcal S}_6}\,   (s_{w_1}-s_{w_2})^2\,
          (s_{w_3}-s_{w_4})^2\,(s_{w_5}-s_{w_6})^2\,,
$$
which vanishes only when four of the $s_i$ coincide.
In particular, this implies that the number of
real solutions does not depend upon the choice of the $s_i$ (when the $s_i$
are distinct). 
Hence both solutions are always real.
When $n=4$ and $\lambda=\hat{1}$, we have checked that, for
each of the  1,001 choices of $s_1,\ldots,s_{10}$ chosen from
$$
  \{ 1,2,3,5,7,10,11,13,15,16,17,23,29,31\}\,,
$$
there are twelve ($=\deg(\hat{1})$) solutions, and all are real.

\section{The Lagrangian Grassmannian}
The {\it Lagrangian Grassmannian} $LG(n)$
is the space of all Lagrangian (maximal isotropic) subspaces in a
$2n$-dimensional vector space $V$ equipped with a nondegenerate alternating
form $\Span{\cdot,\cdot}$.
Such Lagrangian subspaces have dimension $n$.
In contrast to the flag manifolds $\Fla$ and orthogonal Grassmannian
$OG(n)$, we show that there may be no real solutions for the enumerative
problems we consider.
We may assume that $V$ has a ${\mathbb K}$-basis $e_1,\ldots,e_{2n}$, for
which our form is 
$$
  \left\langle \sum x_ie_i,\ \sum y_j e_j\right\rangle\ =\ 
   \sum_{i=1}^n  x_i y_{2n+1-i} - y_i x_{2n+1-i} \,.
$$

An {\it isotropic flag} is a complete flag $\Fdot$ of $V$ such that 
$F_n$ is Lagrangian, and for every $i> n$,
$F_i$ is the annihilator of $F_{2n-i}$; that is,
$\Span{F_{2n-i},\,F_i}\equiv 0$.
An isotropic flag induces a Bruhat decomposition of 
$$
  LG(n)\ =\ \coprod X_\lambda^\circ \Fdot
$$
indexed by strict partitions $\lambda\in SP(n)$.
The Schubert variety $X_\lambda\Fdot$ is the closure of the Schubert cell
$X_\lambda^\circ\Fdot$ and has dimension 
$|\lambda|$.
The Bruhat order is given (as for $OG(n)$) by componentwise comparison of
sequences.  
Although $OG(n)$ and $LG(n)$ have
identical Bruhat decompositions, they are very different spaces.

The unique simple Schubert variety of $LG(n)$ is 
$$
  Y\!\Fdot\ :=\ \{H\in LG(n)\mid H\cap F_n\neq \{0\}\}\,.
$$
Thus $Y\!\Fdot$ is the set-theoretic intersection of $LG(n)$ with the simple
Schubert variety $\Upsilon\Fdot$  of the ordinary Grassmannian $Gr(n)$ of
$n$-dimensional subspaces of $V$.
This is generically transverse.
As with $OG(n)$, the strict partition indexing $Y\!\Fdot$ is
$n,n-1,\ldots,2$. 
The Bruhat decomposition of the Lagrangian Grassmannian is related to 
that of the ordinary Grassmannian in the same way as that of the
orthogonal Grassmannian (see Lemma~\ref{lem:BO-rel}).

Let ${\mathbb K}={\mathbb C}$.
We study families of Schubert varieties defined by isotropic flags
$\Fdot(s)$ osculating a real rational normal curve $\gamma\colon{\mathbb C}\to V$
at $\gamma(s)$.
With our given form $\Span{\cdot,\cdot}$ and basis $e_1,\ldots,e_{2n}$,
one choice for $\gamma$ whose osculating flags are isotropic is
\begin{equation}\label{eq:gamma}
  \gamma(s)\ =\ \left(
      1,\,s,\,\frac{s^2}{2},\,\ldots,\, \frac{s^n}{n!},\,
  -\frac{s^{n+1}}{(n+1)!},\,\frac{s^{n+2}}{(n+2)!},\,
      \ldots,\,(-1)^{n-1}\frac{s^{2n-1}}{(2n-1)!}\right)\,.
\end{equation}

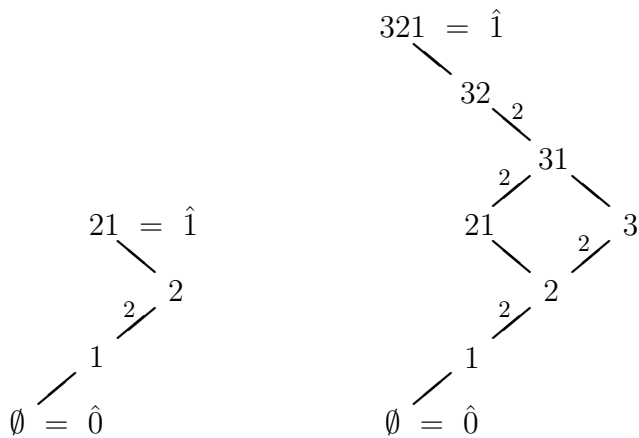
\begin{figure}[htb]
$$
\begin{picture}(95,163)(-1,-3)
\thicklines
\put(55,61){\line(-6,5){14}}
\put(41, 36){\line(6,5){14}}
\put(11, 11){\line(6,5){14}}

\put(43,44){\scriptsize $2$}

\put(30, 75){$21\ =\ \hat{1}$}
\put(60, 50){$2$}
\put(30, 25){$1$}
\put( 0,  0){$\emptyset\ =\ \hat{0}$}
\end{picture}
\qquad\qquad
\begin{picture}(95,163)(-1,-3)
\thicklines
                            \put(25,136){\line(-6,5){14}}
                            \put(55,111){\line(-6,5){14}}
\put(41, 86){\line(6,5){14}}\put(85,86){\line(-6,5){14}}
\put(71, 61){\line(6,5){14}}\put(55,61){\line(-6,5){14}}
\put(41, 36){\line(6,5){14}}
\put(11, 11){\line(6,5){14}}

\put(48,119){\scriptsize $2$}
\put(43, 94){\scriptsize $2$}
\put(73, 69){\scriptsize $2$}
\put(43, 44){\scriptsize $2$}

\put(-2,150){$321\ =\ \hat{1}$}
\put(28.5,125){$32$}
\put(58,100){$31$}
\put(30, 75){$21$}\put(90,75){$3$}
\put(60, 50){$2$}
\put(30, 25){$1$}
\put( 0,  0){$\emptyset\ =\ \hat{0}$}
\end{picture}
$$
 \caption{The multiplicity posets $LG(2)$ and $LG(3)$\label{fig:MP}} 
\end{figure}

Let $\tau$ be the cohomology class dual to the fundamental cycle of
$Y\!\Fdot$, and let $\sigma_\lambda$ be the class dual to the fundamental
cycle of $X_\lambda\Fdot$. 
The Chevalley formula  for $LG(n)$ is
$$
  \sigma_\lambda\cdot \tau\ = \ 
   \sum_{\mu\lessdot\lambda} m^\mu_\lambda\,\sigma_\mu\,,
$$
where the multiplicity $m^\mu_\lambda$ is either 2 or 1, depending
(respectively) upon
whether or not the sequences $\lambda$ and $\mu$ have the same length.
Figure~\ref{fig:MP} shows the multiplicity posets for the enumerative
problem in $LG(2)$ and $LG(3)$ given by the simple Schubert varieties
$Y\!\Fdot(s)$. 

As in Sections~2 and~3, the family ${\mathcal Y}\to{\mathbb A}^1$ whose
fibres are the simple Schubert varieties $Y\!\Fdot(s)$ respects the Bruhat
decomposition of $LG(n)$, and any collection is in general position.
{}From the Chevalley formula, we see that it is not multiplicity-free.

\begin{thm}\label{thm:lagr-families}
The family ${\mathcal Y}\to{\mathbb A}^1$ of simple Schubert
varieties $Y\!\Fdot(s)$ respects the Bruhat
decomposition of $LG(n)$ induced by the flag $\Fdot(0)$.

Any collection of families of simple Schubert varieties is in general
position. 
\end{thm}

The proof of Theorem 4.1, like that of Theorem~\ref{thm:orth-families}, is
virtually identical to that of Theorem~\ref{thm:family}; hence we omit it.

Since the family ${\mathcal Y}$ is not multiplicity-free, we do not have
analogs of Corollaries~\ref{cor:flag-reality} and~\ref{cor:orth-reality}
showing that all solutions may be real.
When $|\lambda|>1$, every chain~(\ref{eq:chains}) in the
multiplicity poset contains the cover $1 < 2$, which has multiplicity 2 and
so is even.
Thus the refined statement of Theorem~\ref{thm:bound} does not
guarantee any real solutions.
We show that there may be no real solutions.

\begin{thm}\label{thm:lagr-nr}
Let $\lambda$ be a strict partition with $|\lambda|=r>1$.
Then there exist real numbers $s_1,\ldots,s_{r}$ such that  the
intersection  
\begin{equation}\label{eq:lagr-int}
  X_\lambda\Fdot(0)\cap Y\!\Fdot(s_1)\cap\cdots\cap
   Y\!\Fdot(s_r)
\end{equation}
is $0$-dimensional and has no real points.
\end{thm}

When $|\lambda|$ is 0 or 1, the degree ($\deg(\lambda)$) of the
intersection~(\ref{eq:lagr-int}) is 1 and so its only point is real.
For all other $\lambda$, $\deg(\lambda)$ is even.
Thus we cannot deduce that the intersection is transverse even for generic
complex numbers $s_1,\ldots,s_{|\lambda|}$.
However, the intersection has been transverse in every case we have computed.
\medskip

\noindent{\bf Proof. }
We induct on the dimension $|\lambda|$ of $X_\lambda\Fdot(0)$ with the
initial case of $|\lambda|=2$ proven in Example~\ref{ex:lagr-nr} (to follow). 
Suppose we have constructed $s_1,\ldots,s_{r-1}\in{\mathbb R}$
having the properties that: (a) for any $\mu$, the intersection
$$
  Y\!\Fdot(s_1)\cap \cdots\cap Y\!\Fdot(s_{r-1})\cap X_\mu\Fdot(0)
$$
is proper; and (b) when $|\mu|=r-1$, it is (necessarily) 0-dimensional, 
has degree $\deg(\mu)$,  and no real points.

Let $\lambda$ be a strict partition with $|\lambda|=r$.
Then the cycle 
$$
   Y\!\Fdot(s_1)\cap \cdots\cap Y\!\Fdot(s_{r-1})\cap
  \sum_{\mu\lessdot \lambda} m^\mu_\lambda\, X_\mu\Fdot(0)
$$
is 0-dimensional, has degree $\deg(\lambda)$, and no real points.
Since the family $Y\!\Fdot(s)$ respects the Bruhat decomposition given by
the flag $\Fdot(0)$, we have
$$
   \lim_{s\to 0} \left( Y\!\Fdot(s) \cap X_\lambda\Fdot(0) \right)
   \ =\ \sum_{\mu\lessdot \lambda} 
    m^\mu_\lambda\, X_\mu\Fdot(0)\,.
$$
Hence there is some $\epsilon_\lambda>0$ such that, if
$0<s_r\leq\epsilon_\lambda$, 
then the  intersection~(\ref{eq:lagr-int}) has dimension 0, degree
$\deg(\lambda)$, and no real points.

Set $s_r=\min\{\epsilon_\lambda: |\lambda|=r\}$.
Since it is an open condition (in the usual topology) on 
$(s_1,\ldots,s_r)\in{\mathbb R}^r$ for the
intersection~(\ref{eq:lagr-int}) to be proper with no real points
and since there are finitely many strict partitions, we may (if necessary) 
choose a nearby $r$-tuple of points such that the
intersection~(\ref{eq:lagr-int}) is proper for every strict partition $\lambda$.  
\qed\medskip
 
\begin{ex}\label{ex:lagr-nr}
When $|\lambda|=2$, we necessarily have $\lambda=2$ and
$$
  X_2\Fdot\ =\ \{H\in LG(n)\mid  F_{n-2}\subset H\subset F_{n+2}\ 
        \mbox{ and }\ \dim (H\cap F_n)\geq n-1\}\,,
$$
which is the image of a simple Schubert variety $Y\Gdot=X_2\Gdot$ of $LG(2)$
under an inclusion $LG(2)\hookrightarrow LG(n)$.
Since $F_{n+2}$ annihilates $F_{n-2}$, the alternating form
$\Span{\cdot,\cdot}$ induces an alternating form on the 4-dimensional space
$W:=F_{n+2}/F_{n-2}$, and the flag $\Fdot$ likewise induces an isotropic
flag $\Gdot$ in $W$.
The inverse image in $F_{n+2}$ of a Lagrangian subspace of $W$ is a
Lagrangian subspace of $V$ contained in $F_{n+2}$.
If we let $\varphi: LG(2)\hookrightarrow LG(n)$ be the induced map, then 
$X_2\Fdot=\varphi(X_2\Gdot)$.

Consider this map for the isotropic flag $\Fdot(\infty)$
of subspaces osculating the point at infinity of  $\gamma$.
Then $(f_1,f_2,f_3,f_4):=(e_{n-1},e_n,e_{n+1},e_{n+2})$ provide a basis for
$W$. 
An explicit calculation using the rational curve $\gamma$~(\ref{eq:gamma})
shows that the flag induced on $W$ is $\Gdot(\infty)$, 
where $\Gdot(s)$ is the flag of subspaces osculating the rational normal
curve $\gamma$ in  $W$ and where 
$\varphi^{-1}(Y\!\Fdot(s))=Y\Gdot(s)$ for $s\in{\mathbb R}$. 
We describe the intersection
$$
  X_2\Fdot(\infty)\cap Y\!\Fdot(s)\cap Y\!\Fdot(t)\ =\ 
  \varphi\left(X_2\Gdot(\infty)\cap Y\Gdot(s)\cap Y\Gdot(t)\right)
$$
when $s$ and $t$ are distinct real numbers.

The Lagrangian subspace $G_2(s)$ is the row space of the matrix
$$
  \left[ \begin{array}{cccc}
        1 & s & s^2/2 & -s^3/6 \\
        0 & 1 &   s   & -s^2/2  
              \end{array}\right]\,.
$$
The flag $\Gdot(\infty)$ is 
$\Span{f_4}\subset\Span{f_4,f_3}\subset\Span{f_4,f_3,f_2}\subset W$.
A Lagrangian subspace in the Schubert cell $X^\circ_2\Gdot(\infty)$
is the row space of the matrix
$$
  \left[ \begin{array}{cccc}
        1 & x &   0   &   y   \\
        0 & 0 &   1   &  -x  
              \end{array}\right]\, ,
$$
where $x$ and $y$ are in ${\mathbb C}$. 
In this way, ${\mathbb C}^2$ gives coordinates for  
the Schubert cell. 
The condition for a Lagrangian subspace $H\in X^\circ_2\Gdot(\infty)$ 
to meet $G_2(s)$, which locally defines the intersection 
$X_2\Gdot(\infty)\cap Y\Gdot(s)$, is
$$
 \det \left[ \begin{array}{cccc}
                  1 & s & s^2/2 & -s^3/6 \\
                  0 & 1 &   s   & -s^2/2 \\
                  1 & x &   0   &   y    \\
                  0 & 0 &   1   &  -x  
                      \end{array}\right]\ \ =\ 
   -y + sx^2 - xs^2 + s^3/3\ =\ 0\,.
$$
If we call this polynomial $g(s)$, then the polynomial
system  $g(s)=g(t)=0$ describes the intersection 
$X_2\Gdot(\infty)\cap Y\Gdot(s)\cap Y\Gdot(t)$.
When $s\neq t$, the solutions are 
$$
  \begin{array}{rcl}
   x&=& {\displaystyle
            \frac{s+t}{2} \pm (s-t)\frac{\sqrt{-3}}{6}}\,,\\
   %
   y &=& {\displaystyle   \rule{0pt}{25pt}
            \frac{s^2t+st^2}{6} \pm (s^2t-st^2)\frac{\sqrt{-3}}{6}}\,,
  \end{array}
$$
which are not real for $s,t\in{\mathbb R}$.

To see that this gives the initial case of Theorem~\ref{thm:lagr-nr}
we observe that, by reparameterizing the rational normal curve, we may move
any three points to any other three points; thus it is no loss to use
$X_2\Fdot(\infty)$ in place of $X_2\Fdot(0)$.
\end{ex}

As before, we ask how much freedom we have to select the
real numbers $s_1,\ldots,s_r$ of Theorem~\ref{thm:lagr-nr}  so
that no points in the intersection~(\ref{eq:orth-int}) are real.
When $n=2$ and $s_1,s_2,s_3$ are distinct and real, no
point in~(\ref{eq:orth-int}) is real.
This is a consequence of Example~\ref{ex:lagr-nr} because, when $n=2$, 
we have $X_2\Fdot=Y\!\Fdot$.
When $n=3$ and $\lambda=\hat{1}$ we have checked that
for each of the  924 choices of $s_1,\ldots,s_6$ chosen from
$$
  \{ 1,2,3,4,5,6,11,12,13,17,19,23\}\,,
$$
there are 16 ($=\deg(\hat{1})$) solutions and none are real.

\section{Schubert Induction for General Schubert Varieties?}
The results in Sections 2, 3, and 4 involve only codimension-1
Schubert varieties
because we cannot show that 
families of general Schubert varieties given by flags
osculating a rational normal curve respect the Bruhat decomposition or that
any collection is in general position. 
Eisenbud and Harris~\cite[Thm.~8.1]{EH83} and~\cite{EH87}
proved this for families $\Omega_\alpha\Fdot(s)$ of arbitrary Schubert
varieties on Grassmannians. 
Their result should extend to all flag manifolds.
We make a precise conjecture for flag varieties of
the classical groups.

Let $V$ be a vector space and $\Span{\cdot,\cdot}$ a bilinear form on $V$,
and set $G:=\mbox{Aut}(V,\Span{\cdot,\cdot})$.
We suppose that $\Span{\cdot,\cdot}$ is either:
\begin{enumerate}
 \item[(1)] identically zero, so that $G$ is a general linear group;
 \item[(2)] nondegenerate and symmetric, so that $G$ is an orthogonal group; or 
 \item[(3)] nondegenerate and alternating, so that $G$ is a symplectic group.
\end{enumerate}
For the orthogonal case, we suppose that $V$ has a basis for
which $\Span{\cdot,\cdot}$ has the  form~(\ref{eq:orth-form}) when $V$
has odd dimension and the same form with $y_{2n+1-i}$ replacing $y_{2n+2-i}$
when $V$ has even dimension.
This last requirement ensures that the real flag manifolds
of $G$ are nonempty.
Let $\gamma$ be a real rational normal curve in $V$ whose flags of
osculating subspaces $\Fdot(s)$ for $s\in\gamma$ are isotropic 
(cases (2) and (3) just listed). 

Let $P$ be a parabolic subgroup of $G$. 
Given a point $0\in\gamma$, the isotropic flag $\Fdot(0)$ induces a
Bruhat decomposition of the flag manifold $G/P$ 
indexed by $w\in W/W_P$, where $W$ is the Weyl group of $G$ and $W_P$ is the
parabolic subgroup associated to $P$. 
For $w\in W/W_P$, let ${\mathcal X}_w\to\gamma$ be the family of Schubert
varieties $X_w\Fdot(s)$.

\begin{conj}\label{conj:Bdec}
For any $w\in W/W_P$, the family ${\mathcal X}_w\to \gamma$ respects the 
Bruhat decomposition of $G/P$ given by the flag
$\Fdot(0)$ and any collection of these families is in general position.
\end{conj}

If this conjecture were true then, for any $u,w\in W/W_P$, we would have
$$
  \lim_{s\rightarrow 0} ( X_u\Fdot(s)\cap X_w )\ =\ 
     \sum_{v\prec w} m^v_{u,\,w}\, X_v\,.
$$
These coefficients $m^v_{u,\,w}$ are the structure constants for the
cohomology ring of $G/P$ with respect to its integral basis of Schubert
classes. 
There are few formulas known for these structure constants, and it is
an open problem to give a combinatorial formula for these
coefficients.
Much of what is known may be found
in~\cite{BS98,BS_lag-pieri,PR96,PR_Pieri_Even,So96}.
An explicit proof of Conjecture~\ref{conj:Bdec} may shed light on
this important problem. 

One class of coefficients for which a formula is known is when $G/P$
is the partial flag manifold $\Fla$ and $u$ is the index of a special
Schubert class.
For these, the coefficient is either 0 or 1~\cite{LS_82,So96}.
A consequence of Conjecture~\ref{conj:Bdec} would be that any
enumerative problem on a partial flag manifold $\Fla$ given by these special
Schubert classes may have all solutions be real, generalizing the result
of~\cite{So99}.

\end{document}